\documentclass[twoside, 12pt]{article}
\usepackage{filecontents}
\usepackage{amsfonts}
\usepackage{amsmath}
\usepackage{amssymb}
\usepackage{blkarray}
\usepackage{multirow}
\usepackage{graphics}
\usepackage[noadjust]{cite}
\usepackage{amsthm} 
\usepackage{varioref}
\usepackage[colorlinks=true,allcolors=blue]{hyperref}
\usepackage[nameinlink,capitalize,noabbrev]{cleveref}
\crefname{equation}{}{}

\usepackage{pgf,tikz}
\usepackage{mathrsfs}
\usetikzlibrary{arrows}
\definecolor{wrwrwr}{rgb}{0.3803921568627451,0.3803921568627451,0.3803921568627451}
\definecolor{rvwvcq}{rgb}{0.08235294117647059,0.396078431372549,0.7529411764705882}
\newtheorem{theorem}{Theorem}[section]
\newtheorem{remark}[theorem]{Remark}
\newtheorem{example}[theorem]{Example}
\newtheorem{lemma}[theorem]{Lemma}
\newtheorem{corollary}[theorem]{Corollary}
\newtheorem{definition}[theorem]{Definition}
\newtheorem{proposition}[theorem]{Proposition}

\newtheorem{conjecture}{Conjecture}

\newcommand{\bea}{\begin{eqnarray}}
\newcommand{\eea}{\end{eqnarray}}

\newcommand{\comment}[1]{}
\def \bpm{\begin{pmatrix}}
\def \epm{\end{pmatrix}}
\def \bd{\begin{definition}}
\def \ed{\end{definition}}
\def \bcc{\begin{conjecture}}
\def \ecc{\end{conjecture}}
\def \bt{\begin{theorem}}
\def \et{\end{theorem}}
\def \bl{\begin{lemma}}
\def \el{\end{lemma}}
\def \bc{\begin{corollary}}
\def \ec{\end{corollary}}
\def\be#1\ee{\begin{align}#1\end{align}}
\def\beq #1\eeq {\begin{align*}#1\end{align*}}
\def \ben{\begin{enumerate}}
\def \een{\end{enumerate}}
\def \ba{\begin{array}}
\def \ea{\end{array}}
\def \bp{\begin{proposition}}
\def \ep{\end{proposition}}
\def \bx{\begin{example}}
\def \ex{\end{example}}
\def \br{\begin{remark}}
\def \er{\end{remark}}
\def \bdsc{\begin{description}}
\def \edsc{\end{description}}
\def\pf{{\it \bf Proof. }}
\def \qed {\hfill \vrule height6pt width6pt depth0pt}
\def\hs{\hspace{.3cm}}

\def\1{1\!\!1}

\begin{document}

\title{Minors of Hermitian (quasi-) Laplacian matrix of a mixed graph}
\author{ Deepak Sarma \\  
Department of Mathematical Sciences,\\ Tezpur University, Tezpur-784028, India.\\
E-mail: \url{deepaks@tezu.ernet.in}
}
\date{}

\pagestyle{myheadings} \markboth{Deepak Sarma}
{Minors of Hermitian (quasi-)Laplacian matrix of a mixed graph}

\maketitle

\vskip 5mm \noindent{\footnotesize

\begin{abstract} A mixed graph is obtained from an unoriented graph by orienting a subset of its edges. Yu, Liu and Qu in 2017 have established the expression for the determinant of the Hermitiann (quasi-) Laplacian matrix of a mixed graph. Here we find general expressions for all minors of Hermitian (quasi-) Laplacian matrix of mixed graphs.

\vskip 3mm

\noindent{\footnotesize Key Words: Mixed Graphs; Hermitian Adjacency Matrix; Hermitian (quasi-) Laplacian matrix }

\vskip 3mm

\noindent {\footnotesize AMS subject classification: 05C20, 05C50, 05B20.}
\end{abstract}

\section{{Introduction}}\label{secone}

\hspace{.62cm}Throughout this article, all graphs are finite and simple, i.e., without loops and multiple edges. But an unoriented edge will be equivalent to two parallel oriented edges in opposite directions. A mixed graph is obtained from an unoriented graph by orienting a subset of its edges. If there is an oriented edge from vertex $u$ to vertex $v$ in a graph $G(V,E)$ then we write it as  $\overrightarrow{uv} \in E$ or $\overleftarrow{vu} \in E$ and we say that $u$ is the head and $v$ is the tail of the edge. If there is an unoriented edge connecting vertices $u$ and $v$ in graph $G(V,E)$ then we write it as $\overline{uv}\in E.$

Hermitian adjacency matrix of a mixed graph was introduced by Liu and Li \cite{ll16} and independently by Guo and Mohar \cite{gm16}.
Let G be a mixed graph with vertex set $V(G)=\{v_{1}, v_{2}, ... , v_{n}\}$ and edge set $E(G)=\{e_{1}, e_{1}, ... , e_{m}\}$. The Hermitian adjacency matrix of the mixed graph $G$ is the $n \times n$ matrix $H(G)=(h_{uv})$, where

\[ h_{uv} = \left\{ \begin{array}{ll}
    i, & \hbox{if $\overrightarrow{uv} \in E$;} \\
    -i, & \hbox{if $\overrightarrow{vu} \in E$;} \\
    1, & \hbox{if $\overline{uv} \in E$;} \\
    0, & \hbox{otherwise.}
  \end{array}
\right. \]

If $V_1\subseteq V(G)$ and $E_1\subseteq E(G),$ then by $G-V_1$ and $G-E_1$ we mean the graphs obtained from $G$ by deleting the vertices in $V_1$ and the edges $E_1$ respectively. In particular case when $V_1=\{u\}$ or $E_1=\{e\},$ we simply write $G-V_1$ by $G-u$ and $G-E_1$ by $G-e$ respectively. If $G$ is an unoriented graph with $u,v\in V(G),$ then the disatance between $u$ and $v,$ denoted by $d_{uv}$ is the length of the shortest path connecting them in $G.$

For mixed graphs, Yu and Qu introduced Hermitian Laplacian matrix \cite{yq15} and Hermitian quasi-Laplacian \cite{ylq17} matrix are denoted and defined as $L_{H}(G)=D(G)-H(G)$ and $Q_{H}(G)=D(G)+H(G)$ respectively,  where $D(G)=diag\{d_{v1}, d_{v2}, \ldots ,d_{vn}\}$ is the diagonal matrix of vertex degrees in the underlying graph. When the graph $G$ is clear from the context we will simply write $L_H$ and $Q_H$ instead of $L_{H}(G)$ and $Q_{H}(G)$ respectively. The authors of \cite{yq15} and \cite{ylq17} have shown that both Hermitian Laplacian and Hermitian quasi-Laplacian matrices of a mixed graph are positive semi-definite. They have also established  expressions for determinants for both Hermitian Laplacian and Hermitian quasi-Laplacian matrices of a mixed graph and gave necessary and sufficient conditions of singularity of the two matrices for any mixed graph.
\par Matrix-tree theorem says that the cofactor of any element of the Laplacian matrix of an unoriented graph equals the number of spanning trees of the graph. In this article we establish the formuli for various minors of Hermitian Laplacian and Hermitian quasi-Laplacian matrices of mixed graphs which generalizes the matrix-tree theorem of unoriented graph. We have organized the paper as that in \cite{bap99}. In section 2, we consider Hermitian incidence matrix and Hermitian quasi-incidence matrix of a mixed graph and find their determinant for
rootless mixed trees and mixed cycles. Also we find expressions for various principal minors of Hermitian Laplacian and Hermitian quasi-Laplacian matrices of mixed graphs. In section 3, we consider non principal submatrices of $Q_H(G)$ and $L_H(G)$ and find their determinants.

\section{Principal Minors}
The authors of \cite{ylq17} introduced two class of incidence matrices of mixed graphs. Here we consider particular cases of those matrices and call them Hermitian incidence matrix and Hermitian quasi-incidences matrix for mixed graphs.

We define the Hermitian quasi-incidence matrix of $G$ as the $n\times m$ matrix $S(G)=(s_{ue})$ whose entries are given by

 \[ s_{ue} = \left\{ \begin{array}{ll}
    1, & \hbox{if $e$ is an unoriented link incident to $u$ or}\\
    & \hbox{if e is an oriented edge with head $u$;} \\
    -i, & \hbox{if e is an oriented edge with tail u;} \\
    0, & \hbox{otherwise.}
  \end{array}
\right. \]

Now to define Hermitian incidence matrix of a mixed graph $G$ we form a new graph $G'$ by assigning arbitrary orientations to unoriented edges in $G$, we call these edges as new oriented edges in $G'$.

Hermitian incidence matrix of $G$ as the $n\times m$ matrix $T(G)=(t_{ue})$ whose entries are given by

 \[ t_{ue} = \left\{ \begin{array}{ll}

    1, & \hbox{if e is an oriented edge in $G$ or new oriented edge in $G'$ with head $u$;} \\
    -1, & \hbox{if e is a new oriented edge in $G'$ with tail $u$;} \\
    i, & \hbox{if e is an oriented edge in $G$ with tail u;} \\
    0, & \hbox{otherwise.}
  \end{array}
\right. \]

When there is no confusion of the graph $G$, we will simply write $S$ in place of $S(G)$ and $T$ in place of $T(G).$

\begin{theorem}\label{t1} \cite{ylq17} For any mixed graph $G,$ $SS^{*}=Q_{H}$ and $TT^{*}=L_{H}$
\end{theorem}

By substructure of a mixed graph we mean an object formed by a subset of the vertex set of the graph together with a subset of the edge set of the graph.
If $G(V,E)$ is a mixed graph and $X \subseteq V$, $Y \subseteq E$, then $R=(X,Y)$ forms a substructure of $G$ and by $S(R)$ and $T(R),$ we denote the corresponding submatrix of quasi-incidence matrix and incidence matrix respectively.
A substructure $R$ of a mixed graph $G$ with equal number of vertices and edges will be called a square substructure. By rootless tree we mean a substructure of a mixed tree obtained by deleting a vertex of the tree. We call the missing vertex as root of the rootless tree and any edge adjacent to the root will be called as rootless edge.  If $T$ is a mixed tree with $v\in V(T)$ and $T_v$ is the rootless tree with root $v$ obtained from $T$, then an edge $e$ in $T_v$ will be treated to be away (towards) the root if $e=\overrightarrow{uw}$ in $T$ and in the underlying tree of $T$ we have $d_{vw}=d_{vu}+1~(d_{vu}=d_{vw}+1). $

\begin{lemma}\label{l1} If the substructure $T_{\diamond}$ of $G$ is a rootless mixed tree, then $det(S(T_{\diamond}))=(-i)^{\alpha},$ where $\alpha$ is the number of directed edges in $T_{\diamond}$ away from the root.
\end{lemma}

\pf First we consider $T_{\diamond}$ to be a rootless mixed path with $n$ vertices and $n$ edges.
Without loss of generality we can assume that edge $e_{n}$ is rootless and its one end is $v_{n}$, each other edge $e_{k}\, (k=1,2,\ldots,n-1)$ is incident with $v_{k}$ and $v_{k+1}$. Then upto permutation similarity the Hermitian quasi-incidence matrix of $T_{\diamond}$ takes the form
$$S(T_{\diamond})=\left(
    \begin{array}{ccccccc}
      * & 0 & 0 & \cdots & 0 & 0 & 0 \\
      * & * & 0 & \cdots & 0 & 0 & 0 \\
      0 & * & * & \cdots & 0 & 0 & 0 \\
      0 & 0 & * & \cdots & 0 & 0 & 0 \\
      \vdots & \vdots & \vdots & \ddots & \vdots & \vdots & \vdots \\
      0 & 0 & 0 & \cdots & * & 0 & 0 \\
      0 & 0 & 0 & \cdots & * & * & 0 \\
      0 & 0 & 0 & \cdots & 0 & * & * \\
    \end{array}
  \right)$$
Here each * is either $1$ or $-i$. In the diagonal, number of $-i$ equals the number of edges away from the root in the path, and so the result follows.

Now we consider $T_{\diamond}$ to be a rootless tree with the edge $e_{n}$ dangling. Let $S(n)$ be the the principal submatrix of $S(T_{\diamond})$ which results after deleting $n^{th}$ row and $n^{th}$ column. Then $S(n)$ represents the Hermitian quasi-incidence matrix of the rootless tree which results from $T_{\diamond}$ after removing $e_{n}$ and the vertex incident to it. So $det(S)=c.det(S(n))$, where $c=-i$ if $e_{n}$ is away from the root and $1$ otherwise. Again applying the same logic and proceeding, ultimately we arrive at some rootless paths, and hence we obtain our required result.

\qed

Similar to the above result we can get the following lemma.

\begin{lemma}\label{l2} If the substructure $T_{\diamond}$ of $G$ is a rootless mixed tree, then $$det(T(T_{\diamond}))=(-1)^ci^{\alpha},$$ where $\alpha$ is the number of directed edges away from the root and c is the number of unoriented edges in $T_{\diamond}.$
\end{lemma}

For a mixed cycle any one direction will be considered as clockwise direction and the opposite direction will be treated as the anticlockwise direction. Now we define five classes of mixed cycles namely type I, type II, type III, type IV and type V so that $\{ type I, type II, type III \}$ and $\{ type I, type IV, type V \}$
forms two partitions of the class of all mixed cycles. Throughout this section we reserve the symbols $a(C)$, $b(C)$ and $c(C)$ respectively for the number of directed edges in clockwise direction, the number of directed edges in anticlockwise direction and the number of unoriented edges in the mixed cycle $C.$ Also when the cycle is understood from the context we would prefer to write them simply as $a,b, c.$

\begin{definition}A mixed cycle will be called a cycle of type I if total number of directed edges is odd.
\end{definition}

\begin{definition}A mixed cycle will be called a cycle of type II if $\frac{|a-b|}{2}+c$ is odd.
\end{definition}

\begin{definition}A mixed cycle will be called a cycle of type III if $\frac{|a-b|}{2}+c$ is even.
\end{definition}

\begin{definition}A mixed cycle will be called a cycle of type IV if $\frac{|a-b|}{2}$ is odd.
\end{definition}

\begin{definition}A mixed cycle will be called a cycle of type V if $\frac{|a-b|}{2}$ is even.
\end{definition}

From the above definitions we can observe the following results for the nature of a mixed cycle when directions of some of its edges are changed.

\begin{theorem}If $C_n$ is a mixed cycle of type III and $C_{n}^{k}$ be the mixed cycle obtained by reverting the directions of k directed edges of $C_n$, then $C_{n}^{k}$ is of type III(type II) if and only if $k$ is even(odd).
\end{theorem}

\begin{theorem} If $C_n$ is a mixed cycle of type V and $C_{n}^{k}$ be the mixed cycle obtained by reverting the directions of k directed edges of $C_n$, then $C_{n}^{k}$ is of type V(type IV) if and only if $k$ is even(odd).
\end{theorem}

\begin{lemma}\label{l3} If $C_n$ is a mixed cycle, then $|det(S(C))|=\sqrt2$ or $2$ or $0$ according as $C$ is of type I or type II or type III respectively

\end{lemma}
\pf Upto permutation similarity, we observe that
the nonzero entries of $S(C_{n})$ occur precisely at positions $(k, k), (k+l, k),$ for $k= 1, 2, \ldots, n-1,$ and at $(1, n)$ and $(n, n).$

$$i.e., \quad S(C_n)=\left(
    \begin{array}{ccccccc}
      * & 0 & 0 & \cdots & 0 & 0 & * \\
      * & * & 0 & \cdots & 0 & 0 & 0 \\
      0 & * & * & \cdots & 0 & 0 & 0 \\
      \vdots & \vdots & \vdots & \ddots & \vdots & \vdots & \vdots \\
      0 & 0 & 0 & \cdots & * & 0 & 0 \\
      0 & 0 & 0 & \cdots & * & * & 0 \\
      0 & 0 & 0 & \cdots & 0 & * & * \\
    \end{array}
  \right)$$

Let $a, b$ respectively denote the number of directed edges in clockwise and anticlockwise direction and $c$ denote the number of unoriented edges in the mixed cycle.
Then expanding along the first row, we see that \beq det(S(C_{n})&=(-i)^{b}+(-1)^{n-1}(-i)^{a}\\
\Rightarrow |det(S(C_{n})|&=|i^{a-b}+(-1)^{n-1}|. \eeq
If $C_{n}$ is of type I, then $a+b$ is odd and so is $a-b.$ Therefore
\beq  i^{a-b}&=\pm i \\
\Rightarrow i^{a-b}+(-1)^{n-1}&=\pm 1\pm i. \eeq

\par Hence $|det(S(C_{n})|= \sqrt{2}  .$ \\
\\
Again if $C_{n}$ is of type II, then $\frac{|a-b|}{2}+c$ is odd. \\
\par Thus $i^{a-b}=(-1)^{c+1}$ and $n$ is even or odd according as $c$ is even or odd. Therefore we get  \\
\beq |det(S(C_{n})|&=|i^{a-b}+(-1)^{c-1}| \\
&= 2\eeq
\\
Finally if $C_{n}$ is of type III, then $\frac{|a-b|}{2}+c$ is even. Thus $n$ and $c$ are either both even or both odd and $i^{a-b}=i^{-2c}.$ Hence

\beq |det(S(C_{n})|=|i^{a-b}+(-1)^{n-1}|
    =|i^{-2c}+(-1)^{c-1}|=0.
\eeq

\qed

\begin{lemma}\label{l4}
If $C_n$ is a mixed cycle, then $|det(T(C))|=\sqrt2$ or $2$ or $0$ according as $C$ is of type I or type IV or type V respectively.
\end{lemma}

\pf Here we observe that if we reverse the direction of any new oriented edge, then the corresponding determinant changes only by sign and thus absolute value remains same. Therefore without loss of generality we consider every new oriented edge in forward direction. Then as in the  \cref{l3}, expanding along the top row, we get
$$|det(T(C_{n})|=|i^{b}+(-1)^{n-1}i^{a}(-1)^c|=|i^{b-a}+(-1)^{n+c-1}|.$$
If the cycle is of type I, we can get our result as that of \cref{l3}. \\

Now if the cycle is not of type I, then $n-c$ is even and therefore $n+c-1$ is odd. Thus we see that\\
 $$ |det(T(C_{n})|=|i^{b-a}-1|.$$
Again if the cycle is of type IV or V then proceeding as in the proof of \cref{l3}, we are done.
\qed

By unicyclic graph of type I/II/III/IV/V, we mean a unicyclic mixed graph where the corresponding cycle is of type I/II/III/IV/V.

\begin{definition} A square substructure $R$ of a mixed graph will be called a special square substructure (in short SSS) if the only possible components of $R$ are rootless trees and/or unicyclic graphs not of type III.
\end{definition}

\begin{lemma}\label{l5} If $R$ is a square substructure of a mixed graph, then

\[ |det(S(R))|=
\left\{
  \begin{array}{ll}
    (\sqrt2)^{x+2y}, & \hbox{if R is an SSS;} \\
   0, & \hbox{otherwise.}
  \end{array}
\right. \]
where $x$ and $y$ denote the number of components of the SSS which are unicyclic graphs of type I and of type II respectively.
\end{lemma}

\pf If R is not an SSS, then it can be observed that every term in the Laplace expansion of det(S(R)) is zero.
\par If R is an SSS, then S(R) is permutationally similar to a block diagonal matrix where each diagonal block corresponds to a component of R. Therefore absolute value of det(S(R)) is the product of the absolute values of the determinants of those diagonal blocks. Also from \cref{l1}, we can see that absolute value of the determinant of quasi-incidence matrix of a rootless tree is always 1. Thus if there are $x$ unicyclic graphs of type I and $y$ unicyclic graphs of type II as a component of R then from \cref{l1} and \cref{l3} we get our required result.
\qed

\begin{definition} A square substructure $R$ of a mixed graph will be called a super special square substructure (in short SSSS) if the only possible components of $R$ are rootless trees and/or unicyclic graphs not of type V.
\end{definition}

\begin{lemma}\label{l6} If $R$ is a square substructure of a mixed graph, then

\[ |det(T(R))|=
\left\{
  \begin{array}{ll}
    (\sqrt2)^{p+2q}, & \hbox{if R is an SSSS;} \\
   0, & \hbox{otherwise.}
  \end{array}
\right. \]
where $p$ and $q$ denote the number of components of the SSSS which are unicyclic graphs of type I and of type IV respectively.
\end{lemma}
\par If $A$ be any $n \times m$ matrix and $B\subseteq \{1, \ldots, m \}$, $C\subseteq \{1, \ldots, n \}$, then $A[B,C]$ will denote the submatrix of $A$ formed by rows corresponding to $B$ and columns corresponding to $C$. Also $A(B,C)$ will denote the submatrix of $A$ formed by deleting rows and columns corresponding to $B$ and $C$ respectively. In short $A[B,B]$ will be denoted by $A[B]$ and $A(B,B)$ will be denoted by $A(B)$.

\begin{theorem}\label{t3} Let $G(V,E)$ be a mixed graph and $W \subseteq V$, then\\
\par det$(Q_{H}[W])=\displaystyle\sum_{R}2^{x+2y}$, where the summation runs over all SSS $R$ with $V(R)=W$ and $x, y$ denote the number of components of R which are unicyclic graphs of type I and of type II respectively.
\end{theorem}

\pf From \cref{t1} we observe that $Q_{H}[W]=S[W,E]S^{*}[W,E]$. Therefore, by Cauchy-Binet Theorem, det$(Q_{H}[W])$ is the sum of modulus of the squares of the determinants of the submatrices S[W, Z], where $Z \subseteq E$ with $|W| = |Z|$. Here we don't need to bother about the zero terms and non zero terms are due to the SSS R corresponding to W. From \cref{l5} we see that every SSS R corresponding to W contributes $((\sqrt2)^{x+2y})^{2}=2^{x+2y}$ to det$(Q_{H}[W])$. Hence the result follows.
\qed

\begin{theorem}\label{lap} Let $G(V,E)$ be a mixed graph and $W \subseteq V$, then\\
\par det$(L_{H}[W])=\displaystyle\sum_{R}2^{p+2q}$, where the summation runs over all SSSS $R$ with $V(R)=W$ and $p, q$ denote the number of components of R which are unicyclic graphs of type I and of type IV respectively.
\end{theorem}

\begin{definition} A mixed graph $G(V,E)$ is said to be quapartite if its vertex set V can be partitioned into four sets $V_{1}$,$V_{2}$,$V_{3}$ and $V_{4}$ so that unoriented edges are between $V_{1}$ and $V_{3}$ or between $V_{2}$ and $V_{4}$ and oriented edges are from $V_{1}$ to $V_{2}$ or $V_{2}$ to $V_{3}$ or $V_{3}$ to $V_{4}$ or $V_{4}$ to $V_{1}$
\end{definition}
\par The above definition generalizes the concept of bipartite graph of unoriented graph to a mixed graph. By close observation of the structure of a quapartite graph, we can easily obtain the following results which extends some well known results of unoriented graph to mixed graph..

\par We know that an unoriented graph is bipartite if and only if all its cycles(if any) are even cycles. We now generalize this concept to mixed graph and provide a pure combinatorial proof for it. If $P=v_0-v_1-\cdots -v_k$ be a mixed path then an edge $\overrightarrow{v_iv_j}\in E(P)$ will be termed as forward or backward in $P$ according as $d_{v_0v_j}=d_{v_0v_i}+1$ or $d_{v_0v_j}=d_{v_0v_i}-1$ in the underlyning graph of $P$. Now we define the weight of a mixed walk as follows.

\begin{definition} If $P$ is a walk in $G$, we define its weight by
$\omega(P)=|\displaystyle\sum_{e\in E(P)}f(e)|,$ where
\[ f(e) = \left\{ \begin{array}{ll}
    1, & \hbox{if $e$ is forward in P;} \\
    -1, & \hbox{if $e$ is backward in P;} \\
    2, & \hbox{if $e$ is unoriented in P;}
  \end{array}
\right. \]
\end{definition}

\begin{lemma}\label{len} If $G$ is a quapartite mixed graph and $u,v$ are in the same vertex partition of $V(G)$, then weight of any walk connecting $u$ and $v$ is an integral multiple of $4$.

\end{lemma}

\pf
If we consider a walk $W$ in a graph $G$ so that except the initial and final vertex of $W$ all other vertices belong to different vertex partition of $V(G)$ then $W$ has length at most 4 and it can be easily observed that $\omega(W)=0 \text{ or }4.$

Now suppose $P=(u=w_1, w_2, \ldots w_i, \ldots w_j, w_{j+1}, \ldots, w_k=v)$ be a $uv$ walk so that $w_i, w_j$ belongs to same vertex partition of $V(G)$ with each $w_s$ belonging to different vertex partition of $V(G)$ for $s=i+1, \ldots, j-1.$ Then by above argument weight of the $w_iw_j$ walk is a multiple of 4. Considering all such possibilities and using above argument we get the required  result.

\qed

\begin{theorem} \label{quaIII} A mixed graph is quapartite if and only if all its cycles are of type III.
\end{theorem}
\pf If $G$ is quapartite and $C$ be any mixed cycle in $G$, then considering $C$ to be a walk connecting any $u\in V(G)$ to itself by \cref{len} we get, $\omega(C)=0(mod 4)$.
Again if $a$, $b$ and $c$ are the number of forward, backward and unoriented edges in $C$, then we get $\omega(C)=a-b+2c.$
Thus $a-b+2c=0(mod 4)$ and hence $C$ is of Type III.

\par Conversely, let every cycle in G are of type III. For any walk P in $G$ we reserve $x$, $y$ and $z$ respectively for the number of forward, backward and unoriented edges in P. Let us call a walk as $ij$ walk if $i$ and $j$ are its terminal vertices. Now for any fix $v\in V(G)$ we consider
\par
$V_{1}=\{u\in V:\frac{|x-y|}{2}+z$ \textit{is even for some uv walk}$\}$
\par $V_{2}=\{u\in V:\frac{|x-y-1|}{2}+z$ \textit{is even for some uv walk}$\}$
\par $V_{3}=\{u\in V:\frac{|x-y|}{2}+z$ \textit{is odd for some uv walk}$\}$
\par $V_{4}=\{u\in V:\frac{|x-y+1|}{2}+z$ \textit{is even for some uv walk}$\}$. \\

\par First we show that $V_{1}$, $V_{2}$, $V_{3}$ and $V_{4}$ actually forms a partition of $V(G)$. \\
Suppose if possible $u\in V_{1}\cap V_{2}$. Then for $u\in V_{1}$ there exists a $uv$ walk $P_{1}$ with $x_{1}$, $y_{1}$ and $z_{1}$ as number of forward, backward and unoriented edges so that
\par
\beq &\frac{|x_{1}-y_{1}|}{2}+z_{1} \text{ is even. }\\
i.e. \quad & i^{x_{1}-y_{1}}=(-1)^{z_{1}}.\eeq
Again for $u\in V_{2}$ there exists a $uv$ walk $P_{2}$ with $x_{2}$, $y_{2}$ and $z_{2}$ as number of forward, backward and unoriented edges so that

\beq &\frac{|x_{2}-y_{2}-1|}{2}+z_{2} \text{ is even. }\\
i.e. \quad & i^{x_{2}-y_{2}-1}=(-1)^{z_{2}}. \eeq

But $P_{1}$ and $P_{2}$ together form a cycle $C$(say). Then for this cycle $a=x_{1}+y_{2}$, $b=y_{1}+x_{2}$ and $c=z_{1}+z_{2}$. Therefore we get

\beq i^{a-b}&= (-1)^{z_{1}}i(-1)^{z_{2}+1}=i(-1)^{c+1} \\
\Rightarrow \hs  i^{a-b-1} &=(-1)^{c+1} \eeq

Which implies that $a-b-1$ is even. Therefore $C$ must be a cycle of type I, a contradiction to our assumption. \\
Similarly we can show disjointness of other pairs of $V_{1}$,$V_{2}$,$V_{3}$ and $V_{4}$.

Again let $u\in V(G)$ be any vertex of $G$ and $P$ be some $uv$ walk in $G$. Then if $|x-y|$ is even then $u\in V_{1}$ or $u\in V_{3}$. Again if $|x-y|$ is odd, then $u\in V_{2}$ or $u\in V_{4}$. Therefore $V_{1}$, $V_{2}$, $V_{3}$ and $V_{4}$ forms a partition of $V(G)$.

\par Now we prove the quapartiteness property of G. For this we consider a fixed edge $uw$ of $G$. Without loss of generality we can take $v\in V(G)$ so that there is a $vu$ walk $P_{1}$ not containing $w.$ Let $P_{2}$ be the walk joining $v$ to $w$ which consist of $P_{1}$ and the edge $uw$. Besides let $x_{1}$, $y_{1}$, $z_{1}$ and $x_{2}$, $y_{2}$, $z_{2}$ in order denote the number of forward, backward and unoriented edges in paths $P_{1}$ and $P_{2}$ respectively.
\par Let us first suppose that $uw$ is an unoriented edge. Then $x_{2}=x_{1}$, $y_{2}=y_{1}$ and $z_{2}=z_{1}+1$.
If $u\in V_{1}$, then for $P_{1}$, $\frac{|x_{1}-y_{1}|}{2}+z_{1}$ is even. \\
Which implies that $\frac{|x_{2}-y_{2}|}{2}+z_{2}$ is odd. Therefore, considering $P_{2}$ we can conclude that $w\in V_{3}$. Similarly if we assume $u\in V_{3}$, then we get $w\in V_{1}$.
\par Now if $u\in V_{2}$, then for $P_{1}$, $\frac{|x_{1}-y_{1}-1|}{2}+z_{1}$ is even. \\
So, $\frac{|x_{2}-y_{2}+1|}{2}+z_{2}=\frac{|x_{1}-y_{1}-1|}{2}+z_{1}+2$ is even.
\par Therefore, considering $P_{2}$ we can conclude that $w\in V_{4}$. \\
Similarly if we assume $u\in V_{4}$, then we get $w\in V_{2}$. \\
Thus we can conclude that unoriented edges are possible only between $V_{1}$ and $V_{3}$ or between $V_{2}$ and $V_{4}$.
\par Next we suppose that $uw$ is an directed edge from vertex $u$ to vertex $w$, Then
\par $x_{2}=x_{1}+1$, $y_{2}=y_{1}$ and $z_{2}=z_{1}$. \\
\\
\textbf{Case 1:} If $u\in V_{1}$, then for $P_{1}$, $\frac{|x_{1}-y_{1}|}{2}+z_{1}$ is even. \\

So, $\frac{|x_{2}-y_{2}-1|}{2}+z_{2}=\frac{|x_{1}-y_{1}|}{2}+z_{1}$ is even.
Therefore $ w\in V_{2}$ \\
\\
\textbf{Case 2:} If $u\in V_{2}$, then for $P_{1}$, $\frac{|x_{1}-y_{1}-1|}{2}+z_{1}$ is even.

Therefore, $\frac{|x_{2}-y_{2}|}{2}+z_{2}=\frac{|x_{1}-y_{1}-1|}{2}+z_{1}+1$, which is odd.
\par Thus considering $P_{2}$ we can conclude that $w\in V_{3}$. \\
\\
\textbf{Case 3:} If $u\in V_{3}$, then for $P_{1}$,
$\frac{|x_{1}-y_{1}|}{2}+z_{1}$ is odd.

Therefore $\frac{|x_{2}-y_{2}+1|}{2}+z_{2}=\frac{|x_{1}-y_{1}|}{2}+z_{1}+1$ is even.
\par Thus considering $P_{2}$ we get $w\in V_{4}$. \\
\\
\textbf{Case 4:} If $u\in V_{4}$, then for $P_{1}$, $\frac{|x_{1}-y_{1}+1|}{2}+z_{1}$ is even. \\
Now $\frac{|x_{2}-y_{2}|}{2}+z_{2}=\frac{|x_{1}-y_{1}+1|}{2}+z_{1}$ is even.
\par Thus considering $P_{2}$ we get $w\in V_{1}$. \\
Hence directed edges are possible only from $V_{1}$ to $V_{2}$ or $V_{2}$ to $V_{3}$ or $V_{3}$ to $V_{4}$ or $V_{4}$ to $V_{1}$. \\
Therefore G is quapartite.
\qed

From the above theorem, we can immediately get the following well known result of graph theory as a corollary.
\begin{corollary}An unoriented graph is bipartite if and only if it does not contain any odd cycle.
\end{corollary}

\begin{theorem} Any odd cycle (cycle with odd number of vertices) in a quapartite graph must consist of both oriented and unoriented edges.
\end{theorem}

\bc
Any directed cycle in a quapartite mixed graph has even number of vertices.
\ec

\bc
Any unoriented cycles in a quapartite mixed graph has even number of vertices.
\ec

\begin{theorem}If G be a quapartite graph and $\widetilde{G}$ be a graph obtained from G by reverting the direction of any one directed edge, then $\widetilde{G}$ is quapartite if and only if that edge does not lie on any cycle.
\end{theorem}
\pf Let $G$ is a quapartite mixed graph and $\widetilde{G}$ be the graph obtained from $G$ by reverting the direction of the directed edge $e\in E(G)$ . Then by \cref{quaIII} all cycles in $G$ are of $type~ III.$ Now if $e$ is not in any cycle of $G,$ then cycles in $\widetilde{G}$ are same as cycles of $G$ and therefore $\widetilde{G}$ is also quapartite. Again if $e$ is in come cycle $C$ in $G$ and $\widetilde{C}$ is the corresponding cycle in $\widetilde{G}$ obtained by reverting the direction of $e,$ then $C$ is of $type ~III$ implies that $\widetilde{G}$ is of $type~II$ in $\widetilde{G}.$ Hence the result follows by \cref{quaIII}. \qed

\begin{theorem}\label{t5} A mixed graph is quapartite if and only if the corresponding Hermitian quasi-Laplacian is singular.
\end{theorem}

\pf If $Q_{H}(G)$ is singular then there exist $x=\{x_{1}, x_{2}, \ldots, x_{n}\}\in C^{n}$ such that $x^{*}S =0.$
Then
\beq
&\overline{x}_{u}+\overline{x}_{v}=0 \text{ for } \overline{uv}\in E \text{ and } \overline{x}_{u}-i\overline{x}_{v}=0 \text{ for } \overrightarrow{uv}\in E. \\
\Rightarrow  \hs & x_{v}=-x_{u} \text{ for } \overline{uv}\in E \text{ and } x_{v}=ix_{u} \text{ for } \overrightarrow{uv}\in E.
\eeq

 \par Now taking,
 \beq
 V_{1}&=\{v_{k}: x_{k}=x_{1}\} \\
 V_{2}&=\{v_{k}: x_{k}=ix_{1}\} \\
 V_{3}&=\{v_{k}: x_{k}=-x_{1}\} \\
 V_{4}&=\{v_{k}: x_{k}=-ix_{1}\} \eeq
 \par we see that $V=V_{1}\cup V_{2}\cup V_{3}\cup V_{4}$ is a partition of $V(G)$ which fulfills the definition of quapartiteness.
 \par Conversely if G is a quapartite mixed graph, then it can be observed that
 $x=(1, 1, \ldots ,1, i, i, \ldots ,i, -1, -1 \ldots ,-1, -i, -i, \ldots ,-i)$ satisfies $x^{*}S=0$. Here $1, i, -1$ and $-i$ appears respectively $|V_{1}|$, $|V_{2}|$, $|V_{3}|$ and $|V_{4}|$ times as components in $x$. Thus $x$ plays the role of an eigenvector for $Q_{H}(G)$ corresponding to the eigenvalue 0. Hence $Q_{H}(G)$ is singular.

\qed

\par The above theorem was proved in \cite{bkp12} in more general form for weighted directed graph, but we have produced the proof here for completeness. Using \cref{t3} and \cref{lap} we can immediately attain the following results.

\bc Quasi-Laplacian matrix of a mixed graph is singular if and only if all its cycles(if any) are of type III.
\ec

\bc Laplacian matrix of a mixed graph is singular if and only if all its cycles(if any) are of type V.
\ec

\section{Non Principal Minors}
Let $G(V,E)$ be any mixed graph. If $Q_{H}[A,B]$ is any non principal square submatrix of $Q_{H}$, then $Q_{H}[A,B]=S(A,E)(S(B,E))^{*}$. Similarly if $L_{H}[A,B]$ is any non principal square submatrix of $L_{H}$, then $L_{H}[A,B]=T(A,E)(T(B,E))^{*}$.  As defined in \cite{bap99}, we would call  $S(A \cup B, F)$ nonsingular relative to A and B if S(A, F) and S(B, F) are both non singular. Similarly we call call  $T(A \cup B, F)$ nonsingular relative to A and B if T(A, F) and T(B, F) are both non singular. If $S(A \cup B, F)$ with $|A|=|B|=|F|$ is non singular then it will be called a quasi generalized matching between A and B.
If $T(A \cup B, F)$ with $|A|=|B|=|F|$ is non singular then it will be called a generalized matching between A and B.

\begin{theorem}\cite{bap99}
Let $G(V, E)$ be a mixed graph; $A, B \subseteq V$ and $F \subseteq E$ with $|A|=|B|=|F|$. Then $S(A \cup B, F)$ is nonsingular relative to $A$ and $B$ if and only if each component is either a nonsingular substructure of $S(A \cap B, F)$ or a tree with exactly one vertex in each of $A \setminus B$ and $B \setminus A$.
\end{theorem}

\begin{theorem}\cite{bap99}
Let $G(V, E)$ be a mixed graph; $A, B \subseteq V$ and $F \subseteq E$ with $|A|=|B|=|F|$. Then $T(A \cup B, F)$ is nonsingular relative to A and B if and only if each component is either a nonsingular substructure of $T(A \cap B, F)$ or a tree with exactly one vertex in each of $A \setminus B$ and $B \setminus A$.
\end{theorem}

\begin{theorem} If $T$ is a tree in a quasi generalized matching R between $A$ and $B$, then it contributes $i^{b-a}$ to the determinant of $Q_{H}[A, B]$, where $a$ and $b$ are the number of edges away from the points of $A\setminus B$ and $B\setminus A$ of $T$ in the path connecting the two points in $T$.
\end{theorem}
\pf Suppose $v_{1}\in A\setminus B$ and $v_{2}\in B\setminus A$ are two points of the tree T. So $T\setminus \{v_{1}\}$ is a rootless tree as a component of an SSS corresponding to $A \subseteq V$ and some $F \subseteq E$ with $|A|=|F|$ and $T\setminus \{v_{2}\}$ is a rootless tree as a component of an SSS corresponding to $B \subseteq V$ and $F \subseteq E$ with $|B|=|F|$. Now $T\setminus \{v_{1}\}$ and $T\setminus \{v_{2}\}$ contributes respectively $(-i)^{\alpha}$ and $(-i)^{\beta}$ to the determinants of S(A,F) and S(B,F), where $\alpha$ and $\beta$ are respectively the number of edges away from $v_{1}$ and $v_{2}$ in the tree T. So they together contribute $(-i)^{\alpha}\times \overline{((-i)^{\beta})}=i^{\beta-\alpha}$ to the determinant of  $S(A\bigcup B, F)$. But if we consider the edges in T which do not lie in the path connecting $v_{1}$ and $v_{2}$, then they all are away from both $v_{1}$ and $v_{2}$, So in our calculation we can avoid them.
Hence a tree T in a quasi generalized matching R between A and B contributes $i^{b-a}$ to the determinant of $Q_{H}[A, B]$, where $a$ and $b$ are the number of edges away from the point of A and B in the path connecting the two points in T.
\qed

Using \cref{l2} and proceeding as in the above theorem, we get the following result.

\begin{theorem} If $T$ is a tree in a generalized matching R between $A$ and $B$, then it contributes $i^{a-b}$ to the determinant of $L_{H}[A, B]$, where $a$ and $b$ are the number of edges away from the points of $A\setminus B$ and $B\setminus A$ of $T$ in the path connecting the two points in $T$.
\end{theorem}

\begin{theorem} If $G(V,E)$ be a mixed graph and $A, B \subseteq V$ with $|A|=|B|$, then $$det(Q_{H}[A, B])={\displaystyle\sum_{R}i^{\Sigma_{T}(b-a)}2^{x+2y}},$$ where the main summation runs over all quasi generalized matching $R$ between A and B and the exponent summation runs over all trees $T$ in $R$, $a$ and $b$ are the number of edges away from the points of $A\setminus B$ and $B\setminus A$ of $T$ in the path connecting the two points in $T$, $x$ and $y$ denote the number of components of $R$ which are unicyclic graphs of type I and of type II respectively.
\end{theorem}

\begin{theorem} If $G(V,E)$ is a mixed graph and $A, B \subseteq V$ with $|A|=|B|$, then $$det(L_{H}[A, B])=\displaystyle\sum_{R}i^{\Sigma_{T}(b-a)}2^{p+2q},$$ where the main summation runs over all generalized matching $R$ between A and B and the exponent summation runs over all trees $T$ in $R$, $a$ and $b$ are the number of edges away from the points of $A\setminus B$ and $B\setminus A$ of $T$ in the path connecting the two points in $T$, $p$ and $q$ denote the number of components of $R$ which are unicyclic graphs of type I and of type IV respectively.
\end{theorem}

\begin{corollary} Let G be a quapartite graph. Then the modulus of all the cofactors of the Hermitian quasi-Laplacian matrix $Q_{H}$ are equal, and their common absolute value is the number of spanning trees of the underlying graph of G.
\end{corollary}

\begin{corollary} Let G be any mixed graph with all cycles(if any) of type V. Then the cofactors of the Hermitian Laplacian matrix $L_{H}$ are equal, and their common value is the number of spanning trees of the underlying graph of G.
\end{corollary}

\begin{corollary}If G is non quapartite mixed graph but $G\setminus v$ is quapartite for some $v\in V(G)$, then $|det(Q_{H}(v))|$ is the number of spanning trees of G.
\end{corollary}

\begin{corollary}If $G(V,E)$ is a quapartite mixed graph and $A, B \subseteq V$ with $|A|=|B|$, then $$det(Q_{H}[A, B])=\displaystyle\sum_{R}i^{\Sigma_{T}(b-a)},$$ where the main summation runs over all generalized matching $R$ between A and B and the exponent summation runs over all trees $T$ in $R$, $a$ and $b$ are the number of edges away from the points of $A\setminus B$ and $B\setminus A$ of $T$ in the path connecting the two points in $T$.
\end{corollary}

\begin{corollary}If $G(V,E)$ be a mixed graph and $A, B \subseteq V$ with $|A|=|B|$ and $A\bigcap B=\phi$, then $$det(Q_{H}[A, B])=\displaystyle\sum_{R}i^{\Sigma_{T}(b-a)},$$ where the main summation runs over all generalized matching $R$ between A and B and the exponent summation runs over all trees $T$ in $R$, $a$ and $b$ are the number of edges away from the points of $A\setminus B$ and $B\setminus A$ of $T$ in the path connecting the two points in $T$.
\end{corollary}
\textbf{Acknowledgement:} The financial assistance for the author was provided by CSIR, India, through JRF.


\end{document}